\documentclass[12pt]{article}

\usepackage{amsmath}
\usepackage{amsfonts}
\usepackage[T1]{fontenc} 
\newtheorem{theo}{Theorem}

\newtheorem{ex}{Example}

\newcommand{\CQFD}{\hfill $\square$}

\def\E{\mathop{\hbox{\rm I\kern-0.20em E}}\nolimits}

\def\og{\leavevmode\raise.3ex
     \hbox{$\scriptscriptstyle\langle\!\langle$~}}
\def\fg{\leavevmode\raise.3ex
     \hbox{~$\!\scriptscriptstyle\,\rangle\!\rangle$}~}
     
\newcommand{\cB}{{\cal B}}

\newcommand{\cF}{{\cal F}}

\newcommand{\cP}{{\cal P}}

\newcommand{\bi}{{\bf i}}
\newfont{\msbm}{msbm10 scaled\magstep1}
\newfont{\msbms}{msbm7 scaled\magstep1} 

\newcommand{\bbD}{\mbox{$\mbox{\msbm D}$}}
\newcommand{\bbE}{\mbox{$\mbox{\msbm E}$}}

\newcommand{\bbP}{\mbox{$\mbox{\msbm P}$}}

\newcommand{\bbR}{\mbox{$\mbox{\msbm R}$}}
\newcommand{\bbS}{\mbox{$\mbox{\msbm S}$}}

\begin{document}

\title{On the convergence of Le Page series in Skohorod space}
\author{Youri Davydov\footnote{Universit\'e des sciences et technologies de Lille, Laboratoire Paul Painlev\'e, UMR CNRS 8524, U.F.R. de Mathématiques, B\^atiment M2, 59655 Villeneuve d'Ascq Cedex, France. Email: Youri.Davydov@math.univ-lille1.fr } \; and Cl\'ement Dombry\footnote{Universit\'e de
Poitiers, Laboratoire LMA, UMR CNRS 6286, T\'el\'eport 2, BP 30179, F-86962 Futuroscope-Chasseneuil cedex,
France.  Email: clement.dombry@math.univ-poitiers.fr}} 
\date{}
\maketitle

\abstract{We consider the problem of the convergence of the so-called Le Page series in the Skohorod space $\bbD^d=\bbD([0,1],\bbR^d)$ and provide a simple criterion based on the moments of the increments of the random process involved in the series. This provides a simple sufficient condition for the existence of an $\alpha$-stable distribution on $\bbD^d$ with given spectral measure.
} 
\\

\section{Introduction}
We are interested in the convergence in the Skohorod space $\bbD^d=\bbD([0,1],\bbR^d)$ endowed with the $J_1$-topology of random series of the form
\begin{equation}\label{eq:LePage}
 X(t)=\sum_{i=1}^\infty \Gamma_i^{-1/\alpha}\varepsilon_i Y_i(t),\quad t\in[0,1],
\end{equation}
where $\alpha\in (0,2)$ and
\begin{itemize}
\item[-] $(\Gamma_i)_{i\geq 1}$ is the increasing enumeration of the points of a Poisson point process on $[0,+\infty)$ with Lebesgue intensity; 
\item[-] $(\varepsilon_i)_{i\geq 1}$ is an i.i.d. sequence of real random variables;
\item[-] $(Y_i)_{i\geq 1}$ is an i.i.d. sequence of $\bbD^d$-valued random variables;
\item[-] the sequences $(\Gamma_i)$, $(\varepsilon_i)$ and $(Y_i)$ are independent.
\end{itemize}
Note that a more constructive definition for the sequence $(\Gamma_i)_{i\geq 1}$ is given by
\[
\Gamma_i=\sum_{j=1}^i \gamma_j,\quad i\geq 1,
\]
where $(\gamma_i)_{i\geq 1}$ is an i.i.d. sequence of random variables with exponential distribution of parameter $1$, and independent of $(\varepsilon_i)$ and $(Y_i)$.

Series of the form \eqref{eq:LePage} are known as Le Page series. For fixed $t\in [0,1]$, the convergence in $\bbR^d$ of the series \eqref{eq:LePage}  is ensured as soon as one of the following conditions is satisfied:
\begin{itemize}
\item[-] $0<\alpha<1$, $\bbE|\varepsilon_1|^\alpha<\infty$ and $\bbE|Y_1(t)|^\alpha<\infty$,
\item[-] $1\leq \alpha <2$, $\bbE\varepsilon_1=0$, $\bbE|\varepsilon_1|^\alpha<\infty$ and $\bbE|Y_1(t)|^\alpha<\infty$.
\end{itemize}
Here $|.|$ denotes the usual Euclidean norm on $\bbR$ or on $\bbR^d$. The random variable $X(t)$ has then an $\alpha$-stable distribution on $\bbR^d$. Conversely, it is well known that any $\alpha$-stable distributions on $\bbR^d$ admits a representation in terms of Le Page series (see for example Samorodnitsky and Taqqu \cite{ST94} section 3.9). 

There is a vast literature on symmetric $\alpha$-stable distributions on separable Banach spaces (see e.g. Ledoux and Talagrand \cite{LT91} or  Araujo and Gin\'e \cite{AG80}). In particular, any symmetric $\alpha$-stable distribution on a separable Banach space can be represented as an almost surely convergent Le Page series (see Corollary 5.5 in \cite{LT91}). The existence of a symmetric $\alpha$-stable distribution with a given spectral measure is not automatic and is linked with the notion of stable type of a Banach space; see Theorem 9.27 in \cite{LT91} for a precise statement. In \cite{DMZ08}, Davydov, Molchanov and Zuyev consider $\alpha$-stable distributions  in the more general framework of abstract convex cones.  

The space $\bbD^d$  equipped with the norm
\[
\|x\|=\sup\{ |x_i(t)|,\quad t\in[0,1],\ i=1,\cdots,d\}, \quad x=(x_1,\cdots,x_d)\in\bbD^d,
\]
is a Banach space but is not separable. The uniform topology associated with this norm is finer than the $J_1$-topology. On the other hand, the space $\bbD^d$ with the $J_1$-topology is Polish, i.e. there exists a metric on $\bbD^d$ compatible with the $J_1$-topology that makes $\bbD^d$ a complete and separable metric space. However, such a metric can not be compatible with the vector space structure since the addition is not continuous in the $J_1$-topology. These properties  explains why the general theory of stable distributions on separable Banach space can not be applied to the space $\bbD^d$. 

Nevertheless, in the case when the series \eqref{eq:LePage} converges, the distribution of the sum $X$ defines an $\alpha$-stable distribution on $\bbD^d$.  We can determine the associated spectral measure $\sigma$ on the unit sphere $\bbS^d=\{x\in\bbD^d;\ \|x\|=1\ \}$. It is given by
\[
\sigma(A)=\frac{\bbE\Big(|\varepsilon_1|^\alpha \|Y_1\|^\alpha \mathbf{1}_{\{\mathrm{sign}(\varepsilon_1)Y_1/\|Y_1\|\in A\}}\Big) } {\bbE(|\varepsilon_1|^\alpha \|Y_1\|^\alpha )},\quad A\in\cB(\bbS^d).
\]
This is closely related to  regular variations theory (see Hult and Lindskog \cite{HL06} or Davis and Mikosch \cite{DM08} ): for all $r>0$ and $A\in\cB(\bbS^d)$ such that $\sigma(\partial A)=0$, it holds that 
\[
\lim_{n\to\infty} n\bbP\Big(\frac{X}{\|X\|}\in A\, \Big|\, \|X\|>rb_n\Big)=r^{-\alpha}\sigma(A),
\]
with 
\[
b_n=\inf\{r>0;\ \bbP(\|X\|<r)\leq n^{-1}\},\quad n\geq 1.
\]
The random variable $X$ is said to be regularly varying in $\bbD^d$ with index $\alpha$ and spectral measure $\sigma$. 

In this framework, convergence of the Le Page series \eqref{eq:LePage} in $\bbD^d$ is known in some particular cases only:
\begin{itemize}
\item[-] When $0<\alpha<1$, $\bbE|\varepsilon_1|^\alpha<\infty$ and $\bbE\|Y_1\|^\alpha<\infty$, the series \eqref{eq:LePage} converges almost surely uniformly in $[0,1]$ (see example 4.2 in Davis and Mikoch \cite{DM08});
\item[-] When $1\leq \alpha <2$, the distribution of the $\varepsilon_i$'s is symmetric, $\bbE|\varepsilon_1|^\alpha<\infty$ and  $Y_i(t)=\mathbf{1}_{[0,t]}(U)$ with $(U_i)_{i\geq 1}$ an i.i.d. sequence of random variables with uniform distribution on $[0,1]$, the series \eqref{eq:LePage} converges almost surely uniformly on $[0,1]$  and the limit process $X$ is a symmetric $\alpha$-stable Lévy process (see Rosinski \cite{R01}).
\end{itemize}

The purpose of this note is to complete these results and to provide a general criterion for almost sure convergence in $\bbD^d$ of the random series \eqref{eq:LePage}. Our main result is the following: 
\begin{theo}\label{theo:main}
Suppose that $1\leq \alpha<2$,  
\[
\bbE\varepsilon_1=0 \quad,\quad\bbE|\varepsilon_1|^\alpha <\infty\quad \mathrm{and}\quad \bbE\|Y_1\|^\alpha<\infty.
\]
Suppose furthermore that there exist $\beta_1,\beta_2>\frac{1}{2}$ and $F_1$ ,$F_2$ nondecreasing continuous functions  on $[0,1]$ such that, for all $0\leq t_1\leq t \leq t_2\leq 1$,
\begin{equation}\label{eq:C1}
\bbE|Y_1(t_2)-Y_1(t_1)|^2\leq |F_1(t_2)-F_1(t_1)|^{\beta_1},
\end{equation}
\begin{equation}\label{eq:C2}
\bbE|Y_1(t_2)-Y_1(t)|^2|Y_1(t)-Y_1(t_1)|^2\leq |F_2(t_2)-F_2(t_1)|^{2\beta_2}.
\end{equation}
Then, the Le Page series \eqref{eq:LePage} converges almost surely in $\bbD^d$.
\end{theo}

The proof of this Theorem is detailled in the next section. We provide hereafter a few cases where Theorem \ref{theo:main} can be applied.
\begin{ex}\label{ex:1}
{\rm The example considered by Davis and Mikosh \cite{DM08} follows easily from Theorem \ref{theo:main}: let $U$ be a random variable with uniform distribution on $[0,1]$ and consider $Y_1(t)=\mathbf{1}_{[0,t]}(U)$, $t\in[0,1]$. Then, for $0\leq t_1\leq t\leq t_2\leq 1$,
\[
 \bbE(Y_1(t_2)-Y_1(t_1))^2=t_2-t_1\quad \mathrm{and}\quad \bbE(Y_1(t_2)-Y_1(t))^2(Y_1(t)-Y_1(t_1))^2=0,
\]
so that conditions \eqref{eq:C1} and \eqref{eq:C2} are satisfied.}
\end{ex}

\begin{ex}\label{ex:2}
{\rm Example \ref{ex:1} can be generalized in the following way: let $p\geq 1$, $(U_i)_{1\leq i\leq p}$  independent random variables on $[0,1]$ and $(R_i)_{1\leq i\leq p}$  random variables on $\bbR^d$. Consider
\[
 Y_1(t)=\sum_{i=1}^p R_i \mathbf{1}_{[0,t]}(U_i).
\]
Assume that for each $i\in\{1,\cdots,p\}$, the cumulative distribution function $F_i$ of $U_i$ is continuous on $[0,1]$. Assume furthermore that there is some $M>0$ such that for all $i\in\{1,\cdots,p\}$
\begin{equation}\label{eq:cond}
 \bbE[R_i^4  \ |\ \cF_U]\leq M  \quad \mathrm{almost\ surely},
\end{equation}
where $\cF_U=\sigma(U_1,\cdots,U_p)$. This is for example the case when the $R_i$'s are uniformly bounded by $M^{1/4}$ or when the $R_i$'s have finite fourth moment and are independent of the $U_i$'s. Simple computations entails that under condition \eqref{eq:cond}, it holds for all $0\leq t_1\leq t\leq t_2\leq 1$,
\[
\bbE(Y_1(t_2)-Y_1(t_1))^2\leq  M^{1/2}p^2 |F(t_2)-F(t_1)|^2
\]
and
\[
\bbE(Y_1(t_2)-Y_1(t))^2(Y_1(t)-Y_1(t_1))^2\leq  Mp^4 |F(t_2)-F(t_1)|^4.
\]
with $F(t)=\sum_{i=1}^p F_i(t)$. So conditions \eqref{eq:C1} and \eqref{eq:C2} are satisfied and Theorem \eqref{theo:main} can be applied in this case.
}
\end{ex}

\begin{ex}\label{ex:3}
{\rm 
A further natural example is the case when $Y_1(t)$ is a Poisson process with intensity $\lambda>0$ on $[0,1]$. Then, for all $0\leq t_1\leq t\leq t_2\leq 1$,
\[
\bbE(Y_1(t_2)-Y_1(t_1))^2=  \lambda|t_2-t_1|+\lambda^2|t_2-t_1|^2
\]
and
\[
\bbE(Y_1(t_2)-Y_1(t))^2(Y_1(t)-Y_1(t_1))^2=(\lambda|t_2-t|+\lambda^2|t_2-t|^2)(\lambda|t-t_1|+\lambda^2|t-t_1|^2)
\]
and we easily see that conditions \eqref{eq:C1} and \eqref{eq:C2} are satisfied.}
\end{ex}


\section{Proof}
For the sake of clarity, we divide the proof of Theorem \ref{theo:main} into five steps.

\vspace{0.5cm} 
\noindent
{\it Step 1.} According to Lemma 1.5.1 in \cite{ST94}, it holds almost surely that for $k$ large enough
\begin{equation}\label{eq:1.1}
|\Gamma_k^{-1/\alpha}-k^{-1/\alpha}|\leq 2\alpha^{-1} k^{-1/\alpha} \sqrt{\frac{\ln\ln k}{k}}.
\end{equation}
This implies the a.s. convergence of the series
\begin{equation}\label{eq:1.2}
\sum_{i=1}^\infty |\Gamma_i^{-1/\alpha}-i^{-1/\alpha}|\, |\varepsilon_i|\, \|Y_i\| <\infty.
\end{equation}
The series \eqref{eq:1.2} has indeed nonnegative terms, and \eqref{eq:1.1} implies that the following conditionnal expectation is finite, 
\[
\bbE\left[\sum_{i=1}^\infty |\Gamma_i^{-1/\alpha}-i^{-1/\alpha}|\, |\varepsilon_i|\, \|Y_i\| \ \Big|\ \cF_\Gamma\right]=\bbE|\varepsilon_1|\,\bbE\|Y_i\|\, \sum_{i=1}^\infty |\Gamma_i^{-1/\alpha}-i^{-1/\alpha}|
\]
where $\cF_\Gamma=\sigma(\Gamma_i,i\geq 1)$. 

This proves that \eqref{eq:1.2} holds true and it is enough to prove the a.s. convergence in $\bbD^d$ of the series 
\begin{equation}\label{eq:1.3}
Z(t)=\sum_{i=1}^\infty i^{-1/\alpha}\varepsilon_i Y_i(t),\quad t\in[0,1],
\end{equation}

\vspace{0.5cm} 
\noindent
{\it Step 2.} Next, consider
\begin{equation}\label{eq:2.1}
\widetilde Z(t)=\sum_{i=1}^\infty i^{-1/\alpha}\tilde \varepsilon_i  Y_i(t),\quad t\in[0,1].
\end{equation}
with
\[
\tilde\varepsilon_i=\varepsilon_i \mathbf{1}_{\{|\varepsilon_i|^\alpha \leq i \}},\quad i\geq 1.
\]
We prove that the series \eqref{eq:1.3} and \eqref{eq:2.1} differ only by a finite number of terms. We have indeed
\[
\sum_{i=1}^\infty \bbP\left(\tilde\varepsilon_i\neq \varepsilon_i \right)= \sum_{i=1}^\infty \bbP\left(|\varepsilon_i|^\alpha > i \right)
\leq  \bbE|\varepsilon_1|^\alpha <\infty 
\]
and the Borel-Cantelli Lemma implies that almost surely $\tilde\varepsilon_i= \varepsilon_i$ for $i$ large enough.
So, both series \eqref{eq:1.3} and \eqref{eq:2.1} have the same nature and it is enough to prove the convergence in $\bbD^d$ of the series \eqref{eq:2.1}.

\vspace{0.5cm} 
\noindent
{\it Step 3.} As a preliminary for step 4, we prove several estimates involving the moments of the random variables $(\tilde \varepsilon_i)_{i\geq 1}$. First, for all $m>\alpha$, 
\begin{equation}\label{eq:3.4}
C(\alpha,m):=\sum_{i=1}^\infty i^{-m/\alpha}\bbE(|\tilde\varepsilon_i|^m) <\infty.
\end{equation}
We have indeed
\begin{eqnarray*}
C(\alpha,m)&=& \sum_{i=1}^\infty i^{-m/\alpha}\bbE(|\varepsilon_i|^m\mathbf{1}_{\{|\varepsilon_i| \leq i^{1/\alpha}\}})\\
&=& \bbE\left(|\varepsilon_1|^m\sum_{i=1}^\infty i^{-m/\alpha}\mathbf{1}_{\{i\geq |\varepsilon_1|^\alpha \}} \right)\\
&\leq & C\bbE(|\varepsilon_1|^m|\varepsilon_1|^{\alpha-m})= C\bbE(|\varepsilon_1|^\alpha)<\infty
\end{eqnarray*}
where the constant  $C=\sup_{x>0}x^{m/\alpha-1} \sum_{i\geq x} i^{-m/\alpha}$
is finite since for $m>\alpha$  
\[
\lim_{x\to\infty}x^{m/\alpha-1}\sum_{i\geq x}^\infty i^{-m/\alpha}= \frac{\alpha}{m-\alpha}.
\]
Similarly, we also have
\begin{equation}\label{eq:3.4bis}
C(\alpha,1):= \sum_{i=1}^\infty i^{-1/\alpha}|\bbE(\tilde\varepsilon_i)| <\infty.
\end{equation}
Indeed, the assumption $\bbE \varepsilon_i =0$ implies
$ \bbE(\tilde\varepsilon_i)= \bbE(\varepsilon_i \mathbf{1}_{\{|\varepsilon_i|^\alpha> i\}})$. Hence,
\begin{eqnarray*}
 \sum_{i=1}^\infty i^{-1/\alpha}|\bbE(\tilde\varepsilon_i)|&\leq& \sum_{i=1}^\infty i^{-1/\alpha}\bbE(|\varepsilon_1| \mathbf{1}_{\{|\varepsilon_1|> i^{1/\alpha}\}})\\
&=& \bbE\Big( |\varepsilon_1| \sum_{i=1}^{[|\varepsilon_1|^\alpha]} i^{-1/\alpha}\Big)\\
&\leq&  \bbE\Big( |\varepsilon_1| C'(|\varepsilon_1|^\alpha)^{1-1/\alpha}\Big)=  C' \bbE |\varepsilon_1|^\alpha <\infty
\end{eqnarray*}
where the constant $C'=\sup_{x>0} x^{1/\alpha-1}\sum_{i=1}^{[x]} i^{-1/\alpha}$ is finite.

\vspace{0.5cm} 
\noindent
{\it Step 4.} For $n\geq 1$, consider the partial sum  
\begin{equation}\label{eq:3.1}
\widetilde Z_n(t)=\sum_{i=1}^n i^{-1/\alpha}\tilde \varepsilon_i  Y_i(t),\quad t\in[0,1].
\end{equation}
We prove that the sequence of processes $(\widetilde Z_n)_{n\geq 1}$ is tight in $\bbD^d$. Following Theorem 3 in Gikhman and Skohorod \cite{GK04} chapter 6 section 3, it is enough to show that there exists $\beta>1/2$ and a non decreasing continuous function $F$ on $[0,1]$ such that
\begin{equation}\label{eq:3.2}
\bbE|\widetilde Z_n(t_2)-\widetilde Z_n(t)|^2|\widetilde Z_n(t)-\widetilde Z_n(t_1)|^2\leq |F(t_2)-F(t_1)|^{2\beta},
\end{equation}
for all $0\leq t_1\leq t\leq t_2\leq 1$. Remark that in Gikhman and Skohorod \cite{GK04}, the result is stated only for  $F(t)\equiv t$. However, the case of a general continuous non decreasing function $F$ follows  easily  from  a simple change of variable.\\

We use the notations $Y(t)=(Y^{p}(t))_{1\leq p\leq d}$, $[\![1,n]\!]=\{1,\cdots,n\}$ and\\ $\bi=(i_1,i_2,i_3,i_4)\in[\![1,n]\!]^4$. We have
\begin{eqnarray}
&&\bbE|\widetilde Z_n(t_2)-\widetilde Z_n(t)|^2|\widetilde Z_n(t)-\widetilde Z_n(t_1)|^2\nonumber\\
&=&\bbE\Big|\sum_{i=1}^n i^{-1/\alpha} \tilde\varepsilon_i (Y_i(t)-Y_i(t_1))\Big|^2\Big|\sum_{j=1}^n j^{-1/\alpha}\tilde\varepsilon_j(Y_j(t_2)-Y_j(t))\Big|^2\nonumber\\
&=&\sum_{1\leq p,q\leq d}\sum_{\bi\in [\![1,n]\!]^4}  (i_1i_2i_3i_4)^{-1/\alpha} \bbE(\tilde\varepsilon_{i_1}\tilde\varepsilon_{i_2}\tilde\varepsilon_{i_3}\tilde\varepsilon_{i_4}) \bbE[(Y_{i_1}^p(t)-Y_{i_1}^p(t_1))\\
&& \qquad (Y_{i_2}^p(t)-Y_{i_2}^p(t_1))(Y_{i_3}^q(t_2)-Y_{i_3}^q(t))(Y_{i_4}^q(t_2)-Y_{i_4}^q(t))] \\
&\leq& d^2\sum_{\bi\in [\![1,n]\!]^4}  (i_1i_2i_3i_4)^{-1/\alpha} |\bbE(\tilde\varepsilon_{i_1}\tilde\varepsilon_{i_2}\tilde\varepsilon_{i_3}\tilde\varepsilon_{i_4})|  D_{\bi}(t,t_1,t_2)
 \label{eq:3.3}
\end{eqnarray}
where 
\[
 D_{\bi}(t,t_1,t_2)=\bbE|Y_{i_1}(t)-Y_{i_1}(t_1)||Y_{i_2}(t)-Y_{i_2}(t_1)||Y_{i_3}(t_2)-Y_{i_3}(t)||Y_{i_4}(t_2)-Y_{i_4}(t)|. 
\]
Consider $\sim_{\bi}$ the equivalence relation on $\{1,\cdots,4\}$ defined by
\[
 j\sim_{\bi} j' \quad \mathrm{if\ and\ only \ if}\quad i_j=i_{j'}. 
\]
Let $\cP$ be the set of all partitions of $\{1,\cdots,4\}$ and  $\tau(\bi)$ be the partition of $\{1,2,3,4\}$ given by the equivalence classes of $\sim_{\bi}$. We introduce these definitions because, since the $Y_i$'s are i.i.d., the term  $D_{\bi}(t,t_1,t_2)$ depends on $\bi$ only through the associated partition $\tau(\bi)$. For example, if $\tau(\bi)=\{1,2,3,4\}$, i.e. if $i_1=i_2=i_3=i_4$, then
\[
D_{\bi}(t,t_1,t_2)=\bbE|Y_{1}(t)-Y_{1}(t_1)|^2|Y_{1}(t_2)-Y_{1}(t)|^2. 
\]
Or if $\tau(\bi)=\{1\}\cup\{2\}\cup\{3\}\cup\{4\}$, i.e. if the indices $i_1,\cdots,i_4$ are pairwise distinct, then
\[
D_{\bi}(t,t_1,t_2)=(\bbE|Y_{1}(t)-Y_{1}(t_1)|\bbE|Y_{1}(t_2)-Y_{1}(t)|)^2.
\]
For $\tau\in\cP$, we denote by $D_\tau(t,t_1,t_2)$ the common value of the terms $D_{\bi}(t,t_1,t_2)$ corresponding to indices $\bi$ such that $\tau(\bi)=\tau$. Define also
\[
 S_{n,\tau}=\sum_{\bi\in\{1,\cdots,n\}^4; \tau(\bi)=\tau}   (i_1i_2i_3i_4)^{-1/\alpha} 
|\bbE(\tilde\varepsilon_{i_1}\tilde\varepsilon_{i_2}\tilde\varepsilon_{i_3}\tilde\varepsilon_{i_4})|.
\]
With these notations, equation \eqref{eq:3.3} can be rewritten as
\begin{equation}
\bbE|\widetilde Z_n(t_2)-\widetilde Z_n(t)|^2|\widetilde Z_n(t)-\widetilde Z_n(t_1)|^2
\leq d^2\sum_{\tau\in\cP} S_{n,\tau} D_\tau(t,t_1,t_2). \label{eq:3.5}
\end{equation}
Under conditions \eqref{eq:C1} and \eqref{eq:C2},  we will prove that for each $\tau\in\cP$, there exist $\beta_\tau>1/2$, a non decreasing continuous function $F_\tau$ on $[0,1]$ and a constant $S_\tau>0$ such that
\begin{equation}
D_\tau(t,t_1,t_2)\leq |F_\tau(t_1)-F_\tau(t_2)|^{2\beta_\tau},\quad 0\leq t_1\leq t\leq t_2, \label{eq:3.6}
\end{equation}
and
\begin{equation}
S_{n,\tau}\leq S_\tau,\quad n\geq 1. \label{eq:3.7}
\end{equation}
Equations \eqref{eq:3.5},\eqref{eq:3.6} and \eqref{eq:3.7} together imply inequality \eqref{eq:3.2} for some suitable choices of $\beta>1/2$ and $F$.

It remains to prove inequalities \eqref{eq:3.6} and \eqref{eq:3.7}.  If $\tau=\{1,2,3,4\}$,
\[
D_\tau(t,t_1,t_2)\leq \bbE|Y_{1}(t)-Y_{1}(t_1)|^2|Y_{1}(t_2)-Y_{1}(t)|^2\leq |F_2(t_2)-F_2(t_1)|^{2\beta_2} 
\]
and
\[
S_n^\tau=\sum_{i=1}^n   i^{-4/\alpha}\bbE\tilde\varepsilon_{i}^4\leq C(\alpha,4).
\]
If $\tau=\{1\}\cup\{2\}\cup\{3\}\cup\{4\}$, Cauchy-Schwartz inequality entails
\[
D_{\tau}(t,t_1,t_2)\leq (\bbE|Y_{1}(t)-Y_{1}(t_1)|\bbE|Y_{1}(t_2)-Y_{1}(t)|)^2 \leq  |F_1(t_2)-F_1(t_1)|^{2\beta_1}
\]
and
\[
 S_n^\tau\leq\sum_{\bi\in\{1,\cdots,n\}^4; \tau(\bi)=\tau}   (i_1i_2i_3i_4)^{-1/\alpha}|\bbE\tilde\varepsilon_{i_1}||\bbE\tilde\varepsilon_{i_2}||\bbE\tilde\varepsilon_{i_3}||\bbE\tilde\varepsilon_{i_4}| \leq C(\alpha,1)^4.
\]
Similarly, for $\tau=\{1,2,3\}\cup\{4\}$, 
\begin{eqnarray*}
D_{\tau}(t,t_1,t_2)&=& \bbE|Y_{1}(t)-Y_{1}(t_1)|^2|Y_{1}(t_2)-Y_{1}(t)| \bbE|Y_{1}(t_2)-Y_{1}(t)|\\
& \leq&  |F_1(t)-F_1(t_1)|^{\beta_1/2}|F_2(t_2)-F_2(t_1)|^{\beta_2}|F_1(t_2)-F_1(t)|^{\beta_1/2}\\
& \leq& |(F_1+F_2)(t_2)-(F_1+F_2)(t_1)|^{\beta_1+\beta_2}
\end{eqnarray*}
and
\[
S_n^\tau\leq\sum_{1\leq i\neq j\leq n}   (i^3j)^{-1/\alpha}\bbE|\tilde\varepsilon_{i}|^3|\bbE\tilde\varepsilon_{j}|  \leq C(\alpha,3)C(\alpha,3).
\]
or for $\tau=\{1,2\}\cup\{3\}\cup\{4\}$,
\begin{eqnarray*}
D_{\tau}(t,t_1,t_2)&=& \bbE|Y_{1}(t)-Y_{1}(t_1)|^2 (\bbE|Y_{1}(t_2)-Y_{1}(t)|)^2\\
& \leq&  |F_1(t)-F_1(t_1)|^{\beta_1}|F_1(t_2)-F_1(t)|^{\beta_1}\\
& \leq& |F_1(t_2)-F_1(t_1)|^{2\beta_1}
\end{eqnarray*}
and
\[
S_n^\tau\leq \sum_{1\leq i\neq j \neq k \leq n}   (i^2jk)^{-1/\alpha}\bbE|\tilde\varepsilon_{i}|^2|\bbE\tilde\varepsilon_{j}||\bbE\tilde\varepsilon_{k}|  \leq C(\alpha,2)C(\alpha,1)^2.
\]
Similar computations can be checked in all remaining cases. The cardinality of $\cP$ is equal to $13$.

\vspace{0.5cm} 
\noindent
{\it Step 5.} We prove Theorem \ref{theo:main}. For each fixed $t\in[0,1]$, Kolmogorov's three-series Theorem implies that $\widetilde Z_n(t)$ converge almost surely as $n\to\infty$. So the finite-dimensional distributions of $(\widetilde Z_n)_{n\geq 1}$ converge. The tightness in $\bbD^d$ of the sequence has already been proved in step 4, so $(\widetilde Z_n)_{n\geq 1}$ weakly convergence in $\bbD^d$ as $n\to\infty$. We then apply Theorem 1 in Kallenberg \cite{K74} and deduce that $\widetilde Z_n$ converges almost surely in $\bbD^d$. In view of step 1 and step 2, this yields the almost sure convergence of the series \eqref{eq:LePage}.\CQFD

\bibliographystyle{plain}
\bibliography{series}
\end{document}